\documentclass[12pt,a4paper]{article}

\usepackage{latexsym}
\usepackage{amsfonts}

\newcommand{\noin}{\noindent}
\newcommand{\diver}{{{\rm{div}}}}

\newtheorem{theorem}{\bf Theorem}[section]
\newtheorem{example}[theorem]{\bf Example}
\newtheorem{lemma}[theorem]{\bf Lemma}

\newtheorem{remark}[theorem]{\bf Remark}

\begin{document}

\title{On Cheng's  Eigenvalue Comparison Theorems}
\author{G. Pacelli Bessa\and
 J. F\'{a}bio Montenegro\footnote{Both authors were partially supported by a CNPq grant.}}
\date{\today}
\maketitle
\begin{abstract}

 \noin We prove Cheng's eigenvalue comparison
 theorems \cite{kn:cheng1}
  for  geodesic balls within the cut locus under weaker geometric hypothesis, Theorems (\ref{theoremCheng1},
 \ref{theoremCheng2}, \ref{theoremCheng3}) and we also show
 that there are certain geometric rigidity in case of equality of the
 eigenvalues. This rigidity becomes isometric rigidity under upper
 sectional curvature bounds or lower Ricci curvature bounds. We
 construct examples of smooth metrics showing that our results are
 true extensions of Cheng´s theorem. We also
 construct a family of complete smooth metrics on $\mathbb{R}^{n}$ non-isometric to the constant
 sectional curvature $\kappa$ metrics of the simply connected space forms $\mathbb{M}(\kappa)$ such that the
 geodesic balls $B_{\mathbb{R}^{n}}(r)$, $B_{\mathbb{M}^{n}(\kappa)}(r)$
   have the same first eigenvalue and the geodesic spheres $\partial B_{\mathbb{R}^{n}}(s)$ and
    $\partial B_{\mathbb{M}^{n}(\kappa)}(s)$, $0<s\leq r$, have the same mean curvatures.
    In the end we construct
   examples of Riemannian manifolds $M$ with arbitrary topology with positive fundamental tone $\lambda^{\ast}>0$
     that generalize Veeravalli's examples,
   \cite{kn:veeravalli}.

\vspace{0.5cm}

  \vspace{.2cm}

\end{abstract}

\section{Introduction} \hspace{.5cm}Let $M$ be
 a complete $n$-dimensional Riemannian manifold and denote by $B_{M}(p,r)$
 the
geodesic ball with center $p$ and radius $r$ and by
$\lambda_{1}(B_{M}(p, r))$  the first Dirichlet eigenvalue of
$B_{M}(p,r)$. Cheng in
 \cite{kn:cheng1}, using a result of Barta \cite{kn:barta},
  proved that if the sectional
curvature of $M$ is bounded above $K_{M}\leq \kappa$ and
$r<\min\{{\rm inj}(p),\pi/\sqrt \kappa\,\}$, ($\pi/\sqrt
\kappa=\infty$ if $\kappa\leq 0$), then $
\lambda_{1}(B_{M}(p,r))\geq
\lambda_{1}(B_{\mathbb{M}(\kappa)}(r))$,  where
$\mathbb{M}(\kappa)$ denote the simply connected space form  of
constant sectional curvature $\kappa$. Cheng also in
\cite{kn:cheng1}
 proved  that if  the Ricci curvature of $M$ is bounded below  $Ric_{M}\geq
 (n-1)\kappa$ then the reverse  inequality $\lambda_{1}(B_{M}(p, r))\leq
\lambda_{1}(B_{\mathbb{M}(\kappa)}(r))$ holds for $r<inj(p)$. In
\cite{kn:cheng2}, choosing a suitable test function for the
Rayleigh quotient, Cheng improved this later inequality proving
that  if $Ric_{M}\geq
 (n-1)\kappa$ then $\lambda_{1}(B_{M}(p,r))\leq
\lambda_{1}(B_{\mathbb{M}(\kappa)}(r))$ for every $r>0$, with
equality holding (for some $r$) if and only if the geodesic balls
$B_{M}(p,r)$ and $B_{\mathbb{M}(\kappa)}(r)$ are
 isometric and $r<inj(p)$. That  raises the questions of whether it is
   possible to prove Cheng's lower eigenvalue inequality
  beyond the cut locus and  show that the geodesic balls are isometric if they have the same first eigenvalue.
   These questions were addressed
  in
 \cite{kn:bessa-montenegro2} and proven to be true, (under upper sectional curvature
 bounds),
 provided  that  the ($n-1$)-Hausdorff measure
  ${\cal H}^{n-1}(Cut (p)\cap B_{M}(p,r))=0$, where $Cut(p)$ is the cut locus of $p$.
   In this paper  we apply  our version of Barta's
  theorem (Theorem \ref{Barta-generalizado})
   to prove an extension of Cheng's
lower and upper eigenvalues inequalities for geodesic balls within
the cut locus (of its center) without sectional or Ricci curvature
bounds. These inequalities have a
 weaker form of geometric rigidity
in the equality  case and we show with family of examples that
this rigidity is all we can expect for. To state our result,
consider $B_{M}(p,r)\subset M$ and
$B_{\mathbb{M}(\kappa)}(r)\subset \mathbb{M}(\kappa)$ geodesic
balls within the cut locus and let $(t,\theta)\in (0,r]\times
\mathbb{S}^{n-1}$ be geodesic coordinates for $B_{M}(p,r)$ and
$B_{\mathbb{M}(\kappa)}(r)$. Let
 $H_{M}(t,\theta)$ and $H_{\mathbb{M}(\kappa)}(t,\theta)=H_{\mathbb{M}(\kappa)}(t)$ be respectively
  the mean curvatures of the distance spheres
$\partial B_{M}(p,t)$  and $\partial B_{\mathbb{M}(\kappa)}(t)$ at
the point $(t,\theta)$ with respect to the unit  vector field
$-\partial/\partial t$. Our first result is the following theorem.
\begin{theorem}\label{theoremCheng1}
If   $H_{M}(s,\theta)\geq H_{\mathbb{M}(\kappa)}(s)$, for all $
s\in (0,r]$ and all $\theta \in \mathbb{S}^{n-1}$ then
\begin{equation}\label{eq1}\lambda_{1}(B_{M}(p,r))\geq
\lambda_{1}(B_{\mathbb{M}(\kappa)}(r)).\end{equation} If
 $H_{M}(s,\theta)\leq H_{\mathbb{M}(\kappa)}(s)$, for all $ s\in (0,r]$ and
 all
$\theta \in \mathbb{S}^{n-1}$ then
\begin{equation}\label{eq2}\lambda_{1}(B_{M}(p,r))\leq
\lambda_{1}(B_{\mathbb{M}(\kappa)}(r)).\end{equation} Equality in
(\ref{eq1}) or  (\ref{eq2}) holds if and only if
$H_{M}(s,\theta)=H_{\mathbb{M}(\kappa)}(s)$, $\forall\, s\in
(0,r]$ and $\forall\,\theta \in \mathbb{S}^{n-1}$.
 \end{theorem}

Observe that the hypotheses of Theorem (\ref{theoremCheng1}) are
implied by an upper  sectional curvature bound $K_{M}\leq \kappa$
and a lower Ricci curvature bound $Ric_{M}\geq (n-1)\kappa$
respectively. On the other hand we construct examples of smooth
metrics on $\mathbb{R}^{n}=[0,\infty)\times \mathbb{S}^{n-1}$ such
that the radial sectional curvatures is bounded below
$K(x)(\partial t,v)>\kappa$ outside a compact set $(x\in
\mathbb{R}^{n}\setminus B_{\mathbb{R}^{n}}(1))$ but
$H_{M}(s,\theta)\geq H_{\mathbb{M}(\kappa)}(s)$, for all $ s\in
(0,\infty)$ and all $\theta \in \mathbb{S}^{n-1}$, see example
(\ref{example1}). This shows that Theorem (\ref{theoremCheng1}) is
a true extension of Cheng´s eigenvalue comparison theorem (within
the cut locus). The rigidity in case of equality of the
eigenvalues, ($H_{M}(s,\theta)=H_{\mathbb{M}(\kappa)}(s)$,
$\forall\, s\in (0,r]$ and $\forall\,\theta \in
\mathbb{S}^{n-1}$), implies that the balls $B_{M}(p,r)$ and
$B_{\mathbb{M}(\kappa)}(r)$ are isometric if we have that
$K_{M}\leq \kappa $ or $Ric_{M}\geq (n-1)\kappa$. Moreover, if the
metric of $B_{M}(p,r)$ is
 expressed in geodesic coordinates by
 $dt^{2}+f^{2}(t)d\theta^{2}$, $f(0)=0$, $f'(0)=1$, $f(t)>0$ for
 $t>0$ then the rigidity (even without curvature bounds) also implies that the
balls $B_{M}(p,r)$ and $B_{\mathbb{M}(\kappa)}(r)$ are
 isometric, see Remark (\ref{R-examples}). This is the case if the  the dimension of $M$ is
 two. On the other hand
we also construct a family of complete smooth metrics $g(\kappa)$
 on $\mathbb{R}^{n}$,
$\kappa < 0$ such that $g(\kappa)$
 is non isometric to the constant sectional
curvature metric of $\mathbb{M}(\kappa)$ but  the geodesic balls
$B_{ g(\kappa)}(r)$,  and $ B_{\mathbb{M}(\kappa)}(r)$ have the
same  first eigenvalue $\lambda_{1}(B_{\mathbb{M}(\kappa)}(r))$
and  their geodesic spheres of same radius
  have the same mean
curvatures, see examples(\ref{example1.1}). These examples show
that the rigidity stated in Theorem (\ref{theoremCheng1}) in
general is all we can expect without curvature bounds. The proof
we present for Theorem (\ref{theoremCheng1}) in fact proves more,
we have few generalizations in
 section \ref{generalizingCheng},  (see Theorems \ref{theoremCheng2}, \ref{theoremCheng3}).
  We also generalize Veeravalli' s examples \cite{kn:veeravalli},
  see Theorem  (\ref{theoremCheng4}).

\section{Preliminaries\label{proofThm1} }
\hspace{.5cm} A powerful tool to obtain lower bounds for the first
Dirichlet eigenvalue of smooth bounded domains in Riemannian
manifolds is the following theorem  proved by J. Barta in
\cite{kn:barta}.
\begin{theorem}[Barta] \label{barta}Let $\Omega\subset M$ be a
domain with compact
 closure and nonempty smooth boundary $\partial \Omega$. Let $\lambda_{1}(\Omega)$
 be the first Dirichlet eigenvalue
 of $\Omega$. Let $f\in
 C^{2}(\Omega)\cap C^{0}(\overline{\Omega})$ with $f>0$ in
 $\Omega$ and $f\vert\partial \Omega=0$. Then
 \begin{equation}\label{eq3}\sup_{\Omega}(-\displaystyle\frac{\triangle f}{f})
 \geq \lambda_{1}(\Omega)\geq
\inf_{\Omega}(-\displaystyle\frac{\triangle f}{f}).
 \end{equation}
 \end{theorem}
 \begin{remark}\label{remark} The first observation is that to prove the lower inequality in (\ref{eq2}) it is
 necessary only to have that $f>0$ in  $\Omega$. A second observation is
 that each of
  the inequalities (\ref{eq3}) is
 strict unless $f$ is a  first eigenfunction of $\Omega$.
 This observation although trivial is essencial in the proof of
  the rigidity statement in Theorem
 (\ref{theoremCheng1}) and its seems to have passed unobserved by Cheng.
 \end{remark}

For  arbitrary open sets $\Omega $, we proved in
\cite{kn:bessa-montenegro2} the following extension of Barta's
Theorem that  gives lower bounds for fundamental tone
$\lambda^{\ast}(\Omega)$. Recall that the fundamental tone
$\lambda^{\ast}(\Omega)$ of an open set $\Omega$ is given by
$$ \lambda^{\ast}(\Omega )= \inf\left\{\frac{\int_{\Omega}\vert
\nabla f\vert^{2}}{\int_{\Omega} f^{2}},\,f\in
 {L^{2}_{1,0}(\Omega ),f\not\equiv 0}\right\},
$$
where $L^{2}_{1,0}(\Omega )$ is the completion  of
$C^{\infty}_{0}(\Omega )$ with respect to the norm $ \Vert\varphi
\Vert_{\Omega}^2=\int_{\Omega}\varphi^{2}+\int_{\Omega}
\vert\nabla \varphi\vert^2. $

\begin{theorem}\label{Barta-generalizado}Let
$\Omega \subset M$ be  an open subset of Riemannian manifold. Then
\begin{equation}\label{eq4}
 \lambda^{\ast}(\Omega)\geq
 \sup_{{\cal X}(\Omega) }\{\inf_{ \Omega} (\diver X-
\vert X \vert^{2})\}, \end{equation}where ${\cal X}(\Omega)$ is
the set of all vector fields $X$ in $\Omega$ such that
$\smallint_{\Omega}\diver (f X)=0$ for all $ f\in
C^{\infty}_{0}(\Omega)$.  If $\Omega$ is a relatively compact open
set  with smooth boundary then
\begin{equation}\label{eq5}\lambda_{1}(\Omega)=
 \sup_{{\cal X}(\Omega ) }\{\inf_{ \Omega} (\diver X-
\vert X \vert^{2})\}. \end{equation}
\end{theorem}
Both results (Barta's Theorem and Theorem
(\ref{Barta-generalizado})) coincides in
 bounded domains with smooth boundaries, but  the vector field aspect
 of this version reveal  the role of the mean curvatures of
 the distance spheres in the comparisons of eigenvalues.
 \subsection{Proof of Theorem \ref{theoremCheng1}}  Let $(t,\theta)\in (0,r]\times
\mathbb{S}^{n-1}$ be geodesic coordinates for $B_{M}(p,r)$ and
$B_{\mathbb{M}(\kappa)}(r)$ and
 $u:B_{\mathbb{M}(\kappa)}(r)\to \mathbb{R}$ be a positive first
 Dirichlet eigenfunction.
 It is well known $u$ is radial function, i.e.
$u(t,\theta)=u(t)$ and $u'(t)\leq 0$. Observe that
$u(t,\theta)=u(t)$ also  defines a smooth function on
 $B_{M}(p,r)$. Now, consider vector
fields $X_{1}$ on $B_{M}(p, r)$ and $X_{2}$ on
$B_{\mathbb{M}(\kappa)}(r)$ given by
\begin{equation}\label{eq6}\begin{array}{lll}X_{1}(t,\theta)&
=&-\displaystyle\frac{u'(t)}{u(t)}\cdot\displaystyle\frac{\partial_{\,1}
}{\partial
t}(t,\theta), \\
& &\\
X_{2}(t,\theta)&
=&-\displaystyle\frac{u'(t)}{u(t)}\cdot\displaystyle\frac{\partial_{\,2}
}{\partial t}(t,\theta).
\end{array}
\end{equation}Here $\displaystyle\frac{\partial_{\,1}
}{\partial t}$ and  $\displaystyle\frac{\partial_{\,2} }{\partial
t}$ are the radial vector fields in $B_{M}(p,r)$ and
$B_{\mathbb{M}(\kappa)}(r)$ respectively. From now on let us write
$B_{M}(r)$ instead $B_{M}(p,r)$ for simplicity of notation.  Now
we have that
\begin{eqnarray}\label{eq7}-\frac{\triangle_{M} u }{u}=\diver_{M} X_{1}-\vert X_{1}\vert^{2} &=&\diver_{M}
X_{1}-\diver_{\mathbb{M}(\kappa)} X_{2} +\vert X_{2}\vert^{2}
-\vert X_{1}\vert^{2}+\diver_{\mathbb{M}(\kappa)}
X_{2}-\vert X_{2}\vert^{2}\nonumber \\
&=&\diver_{M} X_{1}-\diver_{\mathbb{M}(\kappa)} X_{2} -\frac{\triangle_{\mathbb{M}(\kappa)} u }{u}\\
&=&\diver_{M} X_{1}-\diver_{\mathbb{M}(\kappa)} X_{2}+
\lambda_{1}(B_{\mathbb{M}(\kappa)}(r)),\nonumber
\end{eqnarray}since $\diver_{\mathbb{M}(\kappa)} X_{2}-\vert
X_{2}\vert^{2}=-\displaystyle\frac{\triangle_{\mathbb{M}(\kappa)}
u}{u} = \lambda_{1}(B_{\mathbb{M}(\kappa)}(r))$ and $\vert
X_{1}\vert^{2}=\vert X_{2}\vert^{2}$.

\vspace{.2cm}

\noin  By Theorem (\ref{barta})  or (\ref{Barta-generalizado})
  and by identity (\ref{eq7}) we have that
 \begin{equation}\label{eq8}\lambda_{1}(B_{M}(r))\geq\inf_{(t,\theta)}(\diver_{M} X_{1}-
 \vert X_{1}\vert^{2})\geq \inf _{(t,\theta)}[\diver_{M} X_{1}
 -\diver_{\mathbb{M}(\kappa)} X_{2}]+
\lambda_{1}(B_{\mathbb{M}(\kappa)}(r))
\end{equation}Since $ B_{M}(r)$ is a smooth domain we can apply Barta's Theorem  and using identity  (\ref{eq7})
 we  have that
 \begin{equation}\label{eq9}\lambda_{1}(B_{M}(r))\leq \sup_{(t,\theta)}[\diver_{M} X_{1}-
 \vert X_{1}\vert^{2}]\leq \sup _{(t,\theta)}[\diver_{M} X_{1}
 -\diver_{\mathbb{M}(\kappa)} X_{2}]+
\lambda_{1}(B_{\mathbb{M}(\kappa)}(r))
 \end{equation}

 We will associate the difference $\diver_{M}
X_{1}-\diver_{\mathbb{M}(\kappa)} X_{2}$
 to the
mean curvature of the distance spheres through the
 following well known lemma.
\begin{lemma}\label{vector}Let $M\hookrightarrow \overline{M}$ be a smooth
hypersurface. Let $X$ be a smooth vector field on $\overline{M}$.
Then at $x\in M$ we have that
\begin{equation}\label{eq10}\diver_{\overline{M}}X(x)=\diver_{M} X^{t}(x)-\langle X,
\stackrel{\rightarrow}{H}\rangle(x) +\langle
\overline{\nabla}_{\eta}X,\eta \rangle(x),
\end{equation} where $X^{t}$ is the orthogonal projection of $X$
onto the tangent space $T_{x}M$, $\stackrel{\rightarrow}{H}$ is
the mean curvature vector of $M$ at $x$, $\overline{\nabla}$ is
the Levi-Civita connection of $\overline{M}$ and $\eta\in
T_{x}M^{\perp}$.
\end{lemma} Using this lemma we can compute $\diver_{M} X_{1}-\diver_{\mathbb{M}(\kappa)} X_{2}$ at
points  of $B_{M}(r)$ and of $B_{\mathbb{M}(\kappa)}(r)$ with the
same coordinates
 $(t,\theta)$.
 \begin{eqnarray}\label{eq11}\diver_{M} X_{1}-\diver_{\mathbb{M}(\kappa)}
 X_{2}&=&
 -\langle X_{1},
\stackrel{\rightarrow}{H}_{M}\rangle_{M}+\langle X_{2},
\stackrel{\rightarrow}{H}_{\mathbb{M}(\kappa)}\rangle_{\mathbb{M}(\kappa)}\nonumber
\\ & & \nonumber\\
& &+ \,\left\langle \overline{\nabla}^{\,M}_{\partial_{1}/\partial
t}X_{1},\,\frac{\partial_{1}}{\partial t}\right\rangle_{M}
-\,\left\langle
\overline{\nabla}^{\,\mathbb{M}(\kappa)}_{\partial_{2}/\partial
t}X_{2},\,\frac{\partial_{2}}{\partial
t}\right\rangle_{\mathbb{M}(\kappa)}\nonumber\\
&&\nonumber \\
&=&(-u'/u)(H_{M}-H_{\mathbb{M}(\kappa)})+(u'/u)'-(u'/u)'
 \end{eqnarray}Since $$\left\langle \overline{\nabla}^{\,M}_{\partial_{1}/\partial
t}X_{1},\,\frac{\partial_{1}}{\partial t}\right\rangle_{M}\!\!\!
=\left\langle
\overline{\nabla}^{\,\mathbb{M}(\kappa)}_{\partial_{2}/\partial
t}X_{2},\,\frac{\partial_{2}}{\partial
t}\right\rangle_{\mathbb{M}(\kappa)}\!=(u'/u)'$$ and
 $\stackrel{\rightarrow}{H}_{M}=-H_{M}\cdot \partial_{1}/\partial
 t$ and
$\stackrel{\rightarrow}{H}_{\mathbb{M}(\kappa)}=-H_{\mathbb{M}(\kappa)}\cdot\partial_{2}/\partial
t$. Hence
\begin{equation}\label{eq12}\diver_{M}
X_{1}-\diver_{\mathbb{M}(\kappa)}
X_{2}=(-u'/u)(H_{M}-H_{\mathbb{M}(\kappa)}).\end{equation} Now
recall that $(-u'/u)\geq 0$. If
$(H_{M}-H_{\mathbb{M}(\kappa)})\geq 0$ then (\ref{eq8}) and
(\ref{eq12}) implies (\ref{eq1}). Likewise, if
$(H_{M}-H_{\mathbb{M}(\kappa)})\leq 0$ then (\ref{eq9}) and
(\ref{eq12}) implies (\ref{eq2}). To treat the equality case
observe that the proof we presented was nothing but giving  a
suitable positive function $u$ on $B_{M}(r)$ then applying Barta's
Theorem to find the lower bound for $\inf_{B_{M}(r)}
-(\triangle_{M}u/u)\geq\lambda_{1}(B_{\mathbb{M}(\kappa)}(r)) $.
Now, suppose  that
 $\lambda_{1}(B_{M}(r))=\lambda_{1}(B_{\mathbb{M}(\kappa)}(r))$
 then (\ref{eq8}) implies that
$\lambda_{1}(B_{M}(r))=\inf_{(t,\theta)}(\diver_{M} X_{1}-
 \vert X_{1}\vert^{2})$ and  $\inf _{(t,\theta)}[\diver_{M} X_{1}
 -\diver_{\mathbb{M}(\kappa)} X_{2}]=0$.
The Remark (\ref{remark}) says that  the infimum  (supremum) in
(\ref{eq3}) is achieved by a   positive function $f$ if and only
if the function $f$ is an eigenfunction. Thus
$\lambda_{1}(B_{M}(r))=\inf_{(t,\theta)}(\diver_{M} X_{1}-
 \vert X_{1}\vert^{2})$ is saying  that
the function $u:B_{M}(r)\to \mathbb{R}$ is a positive first
 eigenfunction of $B_{M}(r)$, in particular that $\lambda_{1}(B_{M}(r))=\diver_{M} X_{1}-
 \vert X_{1}\vert^{2}$. From
 (\ref{eq7})  we have   that $\diver_{M} X_{1}
 -\diver_{\mathbb{M}(\kappa)}
 X_{2}=\lambda_{1}(B_{M}(r))-\lambda_{1}(B_{\mathbb{M}(\kappa)}(r))=0$. On the other hand,
 $\diver_{M} X_{1}-\diver_{\mathbb{M}(\kappa)}
X_{2}=(-u'/u)(H_{M}-H_{\mathbb{M}(\kappa)})$ and $u'(t)=0$ if and
only if $t=0$. Therefore we have that
$H_{M}(t,\theta)=H_{\mathbb{M}(\kappa)}(t,\theta)$ for all $t>0$
and all $\theta$. The equality in (\ref{eq2}) is treated in the
same way.

 \section{Generalizations of Theorem \ref{theoremCheng1}
 \label{generalizingCheng}}The first generalization we are going
 to consider is the following. Let
 $M$  be a $n$-dimensional
  complete
 Riemannian manifold and let  $B_{M}(r)\subset M$
be a
 geodesic ball within the cut locus. Consider  $\mathbb{R}^{m}=[0,\infty)\times
 \mathbb{S}^{m}$ with metric $ds^{2}=dt^{2}+g^{2}(t)d\xi^{2}$, where $g:[0,\infty)\to \mathbb{R}$
  is a smooth function
satisfying $g(0)=0$,
 $g'(0)=1$, $g(t)>0$ for $t\in (0,\infty)$. Let $B_{\mathbb{R}^{m}}(r)$
 be a geodesic ball of radius $r$.
 Let $(t,\theta)\in (0,r]\times
\mathbb{S}^{n-1}$ be geodesic coordinates for $B_{M}(r)$ and
$(t,\xi)\in (0,r]\times \mathbb{S}^{m-1}$ be geodesic coordinates
for $B_{\mathbb{R}^{m}}(r)$. Let
 $H_{M}(t,\theta)$ and $H_{\mathbb{R}^{m}}(t,\xi)=H_{\mathbb{R}^{m}}(t)$ be respectively
  the mean curvatures of the distance spheres
$\partial B_{M}(t)$  and $\partial B_{\mathbb{R}^{m}}(t)$ at the
points $(t,\theta)$  and $(t,\xi)$ with respect to the unit
vector field $-\partial/\partial t$.
 \begin{theorem}\label{theoremCheng2}
 If
$H_{M}(s,\theta)\geq H_{\mathbb{R}^{m}}(s)=(m-1)(g'/g)(s)$,
$\forall$ $s\in (0,r]$ and $\theta\in \mathbb{S}^{n-1}$ then

\begin{equation}\label{eq13}\lambda_{1}(B_{M}(r))\geq
\lambda_{1}(B_{\mathbb{R}^{m}}(r)).\end{equation}
 If
$H_{M}(s,\theta)\leq H_{\mathbb{R}^{m}}(s)=(m-1)(g'/g)(s)$,
$\forall$ $s\in (0,r]$ and $\theta\in \mathbb{S}^{n-1}$ then

\begin{equation}\label{eq14}\lambda_{1}(B_{M}(r))\leq
\lambda_{1}(B_{\mathbb{R}^{m}}(r)).\end{equation} Equality in
(\ref{eq13}) or (\ref{eq14}) holds if and only if $n=m$ and
$H_{M}(s,\theta)=H_{\mathbb{R}^{n}}(s)$, $\forall\, s\in (0,r]$
and $\theta \in \mathbb{S}^{n-1}$.
 \end{theorem}

A positive first
 eigenfunction $u$ of a geodesic ball $B_{\mathbb{R}^{m}}(r)$
 within the cut locus is radial ($u(t,\xi)=u(t)$) and
 $u'(t)\leq 0$ with $u'(t)=0 \Leftrightarrow t=0$. See a proof of
 that in \cite{kn:chavel}, pages 40-44. Define $v:B_{M}(r)\to
 \mathbb{R}$ by $v(t,\theta)=u(t)$ and take vector fields $X_{1}$ in $B_{M}(r)$
 and $X_{2}$ in $B_{\mathbb{R}^{m}}(r)$ by
 \begin{equation}\label{eq15}\begin{array}{lll}X_{1}(t,\theta)&
=&-\displaystyle\frac{u'(t)}{u(t)}\cdot\displaystyle\frac{\partial_{\,1}
}{\partial
t}(t,\theta), \\
& &\\
X_{2}(t,\xi)&
=&-\displaystyle\frac{u'(t)}{u(t)}\cdot\displaystyle\frac{\partial_{\,2}
}{\partial t}(t,\xi).
\end{array}
\end{equation}
 Proceeding as in the proof of Theorem (\ref{theoremCheng1})
 \begin{eqnarray}-\frac{\triangle_{M}v}{v}(t,\theta)=(\diver_{M}X_{1}-\vert
 X_{1}\vert^{2})(t,\theta)&=&\diver_{M}X_{1}(t,\theta)-\diver_{\mathbb{R}^{m}}X_{2}(t,\xi)\nonumber
 \\ && \nonumber \\
 &
 +&\diver_{\mathbb{R}^{m}}X_{2}(t,\xi)-\vert
 X_{2}\vert^{2}(t,\xi)\nonumber \\
 && \label{eq16} \\
 & +& \vert
 X_{2}\vert^{2}(t,\xi)-\vert
X_{1}\vert^{2}(t,\theta).\nonumber
 \end{eqnarray}
 Since we have that $(\diver_{\mathbb{R}^{m}}X_{2}-\vert
 X_{2}\vert^{2})(t,\xi)=-\displaystyle\frac{\triangle_{\mathbb{R}^{m}}u}{u}=\lambda_{1}(B_{\mathbb{R}^{m}}(r))$,
   $\vert X_{2}\vert^{2}(t,\xi)-\vert X_{1}\vert^{2}(t,\theta)=0$ and $\diver X_{2}(t,\xi)=\diver X_{2}(t)$. Thus we
derive that
 \begin{equation}\label{eq17}\lambda_{1}(B_{M}(r))\geq\inf_{(t,\theta)}(\diver_{M} X_{1}-
 \vert X_{1}\vert^{2})\geq \inf _{(t,\theta)}[\diver_{M} X_{1}
 -\diver_{\mathbb{R}^{m}}X_{2}(t)]+
\lambda_{1}(B_{\mathbb{R}^{m}}(r)).
\end{equation}
Likewise, we can derive

\begin{equation}\label{eq18}\lambda_{1}(B_{M}(r))\leq\sup_{(t,\theta)}(\diver_{M} X_{1}-
 \vert X_{1}\vert^{2})\leq \sup _{(t,\theta)}[\diver_{M} X_{1}
 -\diver_{\mathbb{R}^{m}}X_{2}(t)]+
\lambda_{1}(B_{\mathbb{R}^{m}}(r)).
\end{equation}
 Then applying Lemma (\ref{vector}) we have that $$\diver_{M}
X_{1}(t,\theta)-\diver_{\mathbb{R}^{m}}
X_{2}(t)=-\frac{u'(t)}{u(t)}(H_{M}(t,\theta)-H_{\mathbb{R}^{m}}(t))$$
If $H_{M}(s,\theta)\geq H_{\mathbb{R}^{m}}(s)$, $\forall$ $s\in
(0,r]$ and $\theta\in \mathbb{S}^{n-1}$ then
$\lambda_{1}(B_{M}(r))\geq\lambda_{1}(B_{\mathbb{R}^{m}}(r)).$ On
the other hand if $H_{M}(s,\theta)\leq H_{\mathbb{R}^{m}}(s)$,
$\forall$ $s\in (0,r]$ and $\theta\in \mathbb{S}^{n-1}$ then
$\lambda_{1}(B_{M}(r))\leq\lambda_{1}(B_{\mathbb{R}^{m}}(r)).$ In
case that
$\lambda_{1}(B_{M}(r))=\lambda_{1}(B_{\mathbb{R}^{m}}(r))$ we have
by (\ref{eq17})
 that $\lambda_{1}(B_{M}(r))=\diver_{M} X_{1}-
 \vert X_{1}\vert^{2}$ and $\diver_{M} X_{1}(s,\theta)-\diver_{\mathbb{R}^{m}}
X_{2}(s)= 0$ for all $s\in (0,r]$ and $\theta \in
\mathbb{S}^{n-1}$. Thus by Remark (\ref{remark}) the function $v$
is a positive eigenfunction of $B_{M}(r)$ and
$H_{M}(s,\theta)=H_{\mathbb{R}^{m}}(s)$ for all $s\leq r$,
$\theta\in \mathbb{S}^{n-1}$. To prove that $m=n$ we proceed as
follows. Let $p$ be the center of the ball $B_{M}(r)$. For fixed
$\theta\in \mathbb{S}^{n-1}\subset T_{p}M$, let $\tau_{t}$ denote
parallel translation by $t$ units along the unique minimal
geodesic $\gamma_{\theta}$ satisfying $\gamma_{\theta}(0)=p$ and
$\gamma_{\theta}'(0)=\theta$. For $\eta\in \theta^{\perp}\subset
T_{p}M$ set ${\cal
R}\eta=\tau_{-t}\{R(\gamma_{\theta}'(t),\tau_{t}\eta)\gamma_{\theta}'(t)\}$,
where $R$ is the Riemannian curvature tensor and set ${\cal
A}(t,\theta)$ the path of linear transformations of
$\theta^{\perp}$ satisfying $ {\cal A}''+{\cal R}{\cal A}=0$ with
initial conditions ${\cal A}(0,\theta)=0$, ${\cal
A}'(0,\theta)=I$. The Riemannian metric of $M$ on the geodesic
ball $B_{M}(r)$ is expressed by $ds^{2}(\exp t\theta)=dt^{2}+\vert
{\cal A}(t,\theta)d\theta\vert^{2}$. Set  $\sqrt{G}(t,\theta)=det
{\cal A}(t,\theta)$. The mean curvature $H_{M}(t,\theta)$ of the
geodesic sphere $\partial B_{M}(t)$ at a point $(t,\theta)$ (with
respect to $-\partial/\partial t$) is given by
$\displaystyle\frac{\sqrt{G}'(t,\theta)}{\sqrt{G}(t,\theta)}$.
Moreover for small $t$  we have the  Taylor expansions
$\sqrt{G}(t,\theta)=t^{n-1}(1-t^{2}Ric(\theta,\theta)/6 +
O(t^{3})). $ See \cite{kn:chavel}, pages 316-317. Thus,
\begin{equation}\displaystyle\frac{\sqrt{G}'(t,\theta)}{\sqrt{G}(t,\theta)}=
\frac{(n-1)-(n+1)t^{2}Ric(\theta,\theta)/6
+O(t^{3})}{t(1-t^{2}Ric(\theta,\theta)/6  +O(t^{3}))}
\end{equation}
On the other hand the metric of $B_{\mathbb{R}^{m}}(r)$ is given
by $dt^{2}+g^{2}(t)d\xi^{2}$, where $g(0)=0$, $g'(0)=1$. The mean
curvature  $H_{\mathbb{R}^{m}}(t,\xi)$ of the geodesic sphere
$\partial B_{\mathbb{R}^{m}}(t)$ at a point $(t,\xi)$ is given by
$\displaystyle(m-1)\frac{g'(t)}{g(t)}$. The Taylor expansion of
$g$ is given by $g(t)=t+g''(0)t^{2}/2+O(t^{3})$. Therefore,
\begin{equation}(m-1)\displaystyle\frac{g'(t)}{g(t)}=(m-1)\frac{1+g''(0)t+O(t^{2})}{t(1+g''(0)t/2+O(t^{2}))}
\end{equation}Now, we have that
$H_{M}(t,\theta)=H_{\mathbb{R}^{m}}(t)$ for all $t\in(0,r]$. Then
\begin{equation}
\frac{(n-1)-(n+1)t^{2}Ric(\theta,\theta)/6
+O(t^{3})}{(1-t^{2}Ric(\theta,\theta)/6
+O(t^{3}))}=(m-1)\frac{1+g''(0)t+O(t^{2})}{(1+g''(0)t/2+O(t^{2}))}
\end{equation}Letting $t\to 0$ we have that $n=m$.
\vspace{.2cm}

 Another generalization of Theorem (\ref{theoremCheng1}) is
obtained considering the incomplete cone over an
$(n-1)$-dimensional compact Riemannian manifold  $(N,dh^{2})$. The
incomplete
 cone $C_{f}(N)$ over $N$ is the Riemannian space
 $C(N)=(0,\infty)\times N$ with metric
$ds_{f}^{2}=dt^{2}+f^{2}(t,x)\,dh^{2}$, where $f:[0,\infty)\times
N\to \mathbb{R}$ is a smooth function satisfying $f(0,x)=0, \,f'
(0,x)=1 $, $f(t,x)>0$ for all $t>0$. The completed cone
$\overline{C_{f}(N)}=C_{f}(N)\cup \{p\}$, $p=\{0\}\times N$. The
Euclidean space  $\mathbb{R}^{m}$ with metric
$ds^{2}=dt^{2}+g^{2}(t)d\theta^{2}$ is the completed cone
$\overline{C_{g}(\mathbb{S}^{m-1})}$. The next theorem compares
the fundamental tone $\lambda^{\ast}(C_{f}(N)(r))$ of the the
trunked cone $C_{f}(N)(r)=(0,r)\times N$ with the lowest Dirichlet
eigenvalue $\lambda_{1}(B_{\mathbb{R}^{m}}(r))$ of the geodesic
ball $B_{\mathbb{R}^{m}}(r)$.

\begin{theorem}\label{theoremCheng3}  Let $C_{f}(N)$ be  a incomplete cone
 over a
compact $(n-1)$-dimensional
 Riemannian  manifold $(N, dh^{2})$ and $\mathbb{R}^{m}$ with
 metric $ds^{2}=dt^{2}+g^{2}(t)d\theta^{2}$. If
 \begin{equation}\label{eq22}(n-1)(f'/f)(t,x)\geq (m-1)(g'/g)(t),
\end{equation}for all $ x\in N$ and all $t\in (0,r)$ where   $'$ means the
derivative with respect to the variable $t$.
 Then
\begin{equation}\label{eq23}\lambda^{\ast}(C_{f}(N)(r)) \geq
\lambda_{1}(B_{\mathbb{R}^{n}}(r))\end{equation} If (\ref{eq22})
holds for all $t>0$ then letting $r\to \infty$ we have that
$$\lambda^{\ast}(C_{f}(N))\geq
\lambda^{\ast}(\mathbb{R}^{m})
$$
\end{theorem} The proof of Theorem  (\ref{theoremCheng3}) is
similar
 to the proof of Theorem (\ref{theoremCheng1}). We take $u$ to be
 a positive first Dirichlet eigenfunction of
 $B_{\mathbb{R}^{n}}(r)$ and consider the vector fields
 $X_{1}(t,x)=-\displaystyle\frac{u'}{u}(t)\cdot
  \displaystyle\frac{\partial_{\,1}}{\partial t}(t,x)$ and
 $X_{2}(t,\theta)=-\displaystyle\frac{u'}{u}(t)\cdot
  \displaystyle\frac{\partial_{\,2}}{\partial t}(t,\theta)$. Thus we
  have by Theorem (\ref{Barta-generalizado}) that
  \begin{equation}\lambda^{\ast}(C_{f}(N)(r)) \geq
 \inf[ \diver X_{1}-\vert X_{1} \vert^{2}]=\inf[\diver X_{1}-\diver X_{2}]+
 \lambda_{1}(B_{\mathbb{R}^{n}}(r)).\label{eq24} \end{equation}
 Observe that the slice ${t}\times N$, $\{t\}\in (0,r)$ is
 a smooth hypersurface of $C_{f}(N)(r)$ thus we may apply Lemma
 (\ref{vector}) to obtain that $\diver X_{1}(t,x)-\diver
 X_{2}(t,\theta)=(n-1)\displaystyle
 \frac{f'}{f}(t,x)-(m-1)\displaystyle\frac{g'}{g}(t)\geq 0$. This together
 with (\ref{eq24}) proves (\ref{eq23}).

 \vspace{.2cm}
\noin These ideas used in the proofs of theorems
(\ref{theoremCheng1}, \ref{theoremCheng2}, \ref{theoremCheng3})
can be used to obtain examples of Riemannian manifolds $M$ with
arbitrary fundamental groups and variable sectional curvatures and
with positive fundamental tone $\lambda^{\ast}(M)>0$. For
instance,  let $M=\mathbb{R}^{m}\times N$ with the metric
$ds^{2}=dt^{2}+f^{2}(t,\theta)d\theta^{2}+g^{2}(t,\theta)dh^{2}$
where $(N,dh^{2})$ is a complete $n$-dimensional Riemannian
manifold and  $f, g:\mathbb{R}^{n}\to [0,\infty)$ are  smooth
functions, $f$ satisfying  $f(0,\theta)=0, \,f' (0,\theta)=1 $ and
$f(t,\theta)>0$ for  $t>0$ and  $\theta\in \mathbb{S}^{m-1}$,
$g(t,\theta)>0$ for all $(t,\theta)\in [0,\infty)\times
\mathbb{S}^{m-1}$. Let $\Omega=B_{\mathbb{R}^{m}}(r)\times
W\subset M$ where $B_{\mathbb{R}^{m}}(r)\subset\mathbb{R}^{m}$ is
a  ball with radius $r$ nd $W\subset N$ is a domain with compact
closure and smooth boundary $\partial W$ (possibly empty). Let
$B_{\mathbb{M}^{\,l}(\kappa)}(r)\subset\mathbb{M}^{\,l}(\kappa)$
be a geodesic ball of radius $r$ in the simply connected
$l$-dimensional space form of constant sectional curvature
$\kappa$ with metric $dt^{2}+S_{\kappa}^{\,2}(t)d\theta^{2}$.
\begin{theorem}\label{theoremCheng4}
 If $
(m-1)\displaystyle\frac{f'}{f}(t,\theta)+n\displaystyle\frac{g'}{g}(t,\theta)\geq
(l-1)\displaystyle\frac{S_{\kappa}'}{S_{\kappa}}(t)$  for all
$t\in [0,r]$ and  $\theta \in \mathbb{S}^{m-1}$,  then
\begin{equation}\label{eq27}\lambda_{1}(\Omega)\geq
\lambda_{1}(B_{\mathbb{M}^{\,l}(\kappa)}(r))+\inf_{(t,\,\theta)\in\,
\Omega}[
\displaystyle\frac{1}{g^{2}}]\cdot\lambda_{1}(W).\end{equation}
  If $r=\infty$ and letting $W=N$ we have that $$ \lambda^{\ast}(M)\geq
 (l-1)^{2}\kappa^{2}/4 +
\inf_{t,\,\theta}[
\displaystyle\frac{1}{g^{2}}]\cdot\lambda^{\ast}(N).$$ If  $
(m-1)\displaystyle\frac{f'}{f}(t,\theta)+n\displaystyle\frac{g'}{g}(t,\theta)\leq
(l-1)\displaystyle\frac{S_{\kappa}'}{S_{\kappa}}(t)$  for all
$t\in [0,r]$ and  $\theta \in \mathbb{S}^{m-1}$,  then

\begin{equation}\label{eq28}\lambda_{1}(\Omega)\leq
\lambda_{1}(B_{\mathbb{M}^{\,l}(\kappa)}(r))+\sup_{(t,\,\theta)\in\,
\Omega}[
\displaystyle\frac{1}{g^{2}}]\cdot\lambda_{1}(W).\end{equation} If
$r=\infty$, and letting $W=N$ we have that
$$\lambda^{\ast}(M)\leq (l-1)^{2}\kappa^{2}/4+ \sup_{(t,\theta)}[
\displaystyle\frac{1}{g^{2}}]\cdot\lambda^{\ast}(N)$$

\end{theorem}

\noin Choose a positive function $\psi:\Omega\to \mathbb{R}$ given
by $\psi(t,\theta,x)=u(t)\cdot \xi(x)$ where  $u$ and $\xi$ are
positive eigenfunctions of $B_{\mathbb{M}^{l}(\kappa)}(r)$ and $W$
respectively, i.e. $u$ satisfies the differential equation
\begin{equation}\label{eq29}\frac{\partial^{\,2}u}{\partial
\,t^{2}}(t)+
(l-1)\displaystyle\frac{S_{k}'}{S_{k}}(t)\frac{\partial
u}{\partial
\,t}(t)+\lambda_{1}(B_{\mathbb{M}^{l}(\kappa)}(r))\,u(t)=0
\end{equation}with $u(0)=1$, $u'(0)=0$ and $\xi:W\to \mathbb{R}$
satisfies $\triangle_{dh^{2}}\xi +\lambda_{1}(W)\xi=0$ in $W$ and
$\xi\vert
\partial W=0$. It is clear that $\psi\in C^{2}(\Omega)\cap
C^{0}(\overline{\Omega_{2}})$ with $\psi>0$ in $\Omega$ and
$\psi\vert
\partial \Omega=0$. The Laplace operator of $ds^{2}$ is written in
geodesic coordinates is given
 by
\begin{eqnarray}\label{eq30}
\triangle_{ds^{2}} & = & \frac{\partial^{2}}{\partial
\,t^{2}}+\left[(m-1)\,\frac{1}{f}\,\frac{\partial f}{\partial t}+
n\,\frac{1} {g}\,\frac{\partial g}{\partial
t}\right]\frac{\partial}{\partial t}+
\frac{m-3}{f^{3}}\,\,d\theta^{2}\left(
\nabla_{d\theta^{2}}\,f,\,\nabla_{d\theta^{2}}\,\cdot\,\right)\nonumber \\
& +& \frac{n}{g f^2  }\,\,d\theta^{2}\left(
\nabla_{d\theta^{2}}\,g,\,\nabla_{d\theta^{2}}\,\cdot\,\right) +
\frac{1}{f^{2}}\triangle_{d\theta^{2}}+
\frac{1}{g^{2}}\triangle_{dh^{2}}
\end{eqnarray} where $ \nabla_{d\theta^{2}}$ and $\triangle_{d\theta^{2}}$  are
 respectively the
gradients and the Laplacian of  $\mathbb{S}^{m-1}$ and $
\nabla_{dh^{2}}$ and $\triangle_{dh^{2}}$  are respectively the
gradients and the Laplacian of $N$. Computing $-\triangle_{ds^{2}}
\psi/ \psi$ we have,
\begin{eqnarray}-\frac{\triangle_{ds^{2}} \psi}{\psi}& =  &
-\frac{u''}{u}+(m-1)\,\frac{f'}{f}\frac{ u'}{u}+ n\,\frac{g'}
{g}\frac{ u'}{u}
-\frac{1}{g^{2}}\frac{\triangle_{dh^{2}}\xi}{\xi}\nonumber \\
&=& \lambda_{1}(B_{\mathbb{M}^{s}(\kappa)}(r)) -\frac{u'}{u}\left(
(m-1)\,\frac{f'}{f}+ n\,\frac{g'}
{g}-(s-1)\frac{S_{k}'}{S_{k}}\right)+\,
\frac{1}{g^{2}}\lambda_{1}(W)\label{eq31}
\end{eqnarray}
 If $\displaystyle(m-1)\,\frac{f'}{f}+ n\,\frac{g'}
{g}\geq (s-1)\frac{S_{k}'}{S_{k}}$  then from (\ref{eq31})
$$\inf(-\frac{\triangle \psi}{\psi})\geq
\lambda_{1}(B_{\mathbb{M}^{s}(\kappa)}(r)) + \inf
\frac{1}{g^{2}}\lambda_{1}(W).
$$

\noin If $\displaystyle(m-1)\,\frac{f'}{f}+ n\,\frac{g'}
{g}-(s-1)\frac{S_{k}'}{S_{k}}\leq 0$ then $$\sup(-\frac{\triangle
\psi}{\psi})\leq \lambda_{1}(B_{\mathbb{M}^{s}(\kappa)}(r)) + \sup
\frac{1}{g^{2}}\lambda_{1}(W). $$Since $-(u'/u)\geq 0$.

\section{Examples \label{examples}}In this section we construct
examples of metrics showing certain aspects  of Cheng's eigenvalue
comparison theorem. In this first example we construct a family of
metrics on $\mathbb{R}^{n}$ with radial sectional  $
K_{\mathbb{R}^{n}}>\kappa$  outside a compact set and such that
the
  mean curvatures of the distance spheres satisfy
$H_{\mathbb{R}^{n}}(t,\theta) \geq
H_{\mathbb{M}^{n}(\kappa)}(t,\theta)=(n-1)(S'_{\kappa}/S_{\kappa})(t)$.
\begin{example}\label{example1}
Let $\mathbb{R}^{n}=[0,\infty)\times \mathbb{S}^{n-1}$ with the
metric $ds^{2}=dt^{2}+f^{2}(t)d\theta^{2}$, $f(0)=0$, $f'(0)=1$.
Set $\psi_{\kappa}(t)=(-f'S_{\kappa}+ fS'_{\kappa})(t)$, where $'$
means differentiation with respect to $t$ and $S_{\kappa}$ is
given by
\begin{equation}\label{eqSk}\begin{array}{ll}S_{\kappa}(t)=\left\{\begin{array}{lll}\sinh
(\sqrt-\kappa
t)/\sqrt-\kappa& if & \kappa =-k^{2}\\
t &if & \kappa=0\\
\sin (\sqrt\kappa t)/\sqrt \kappa& if & \kappa =k^{2}
\end{array}\right.,& C_{\kappa}(t)=S_{\kappa}'(t)\end{array}
\end{equation}
The radial sectional curvature of $(\mathbb{R}^{n}, ds^{2})$ is
bounded above by $\kappa$ if and only  if $\psi_{\kappa}'(t)\leq
0$. The mean curvatures of  $\partial B_{\mathbb{R}^{n}}(t)$  and
$\partial B_{\mathbb{M}^{n}(\kappa)}(t)$ satisfies
 $H_{\mathbb{R}^{n}}(t,\theta)\geq
H_{\mathbb{M}^{n}(\kappa)}(t)$ if and only if
$\psi_{\kappa}(t)\leq 0$. From $\psi_{\kappa}(t)=(-f'S_{\kappa}+
fS'_{\kappa})(t)$ we have that
$\psi_{\kappa}(0)=\psi_{\kappa}'(0)=0$. Solving the differential
equation we have
$$f(t)=S_{\kappa}(t)+S_{\kappa}(t)\int_{0}^{t}\psi_{\kappa}(s)/S_{\kappa}(s)ds$$
Let  $\psi_{\kappa}:[0,\infty)\to \mathbb{R}$ be a  smooth
function satisfying $\psi_{\kappa}(0)=\psi_{\kappa}'(0)=0$, $\psi
(t)\leq 0$,  $\psi_{\kappa}'(t)>0$ for $t>1$ and $\vert
\int_{0}^{t}\psi_{\kappa}(s)/S_{\kappa}(s)ds \vert<\infty$. This
yields a metric $ds_{f}^{2}=dt^{2}+f^{2}(t)d\theta$ with sectional
curvature
 $ K_{\mathbb{R}^{n}}>\kappa$  outside a compact set and
such that the
  mean curvatures of the distance spheres satisfy
$(n-1)(f'/f)(t)=H_{\mathbb{R}^{n}}(t,\theta) \geq
H_{(\mathbb{M}^{n}(\kappa)}(t,\theta)=(n-1)(S'_{\kappa}/S_{\kappa})(t)$.
\end{example}
\begin{remark}\label{R-examples} If the metric of $M$ is expressed by
$dt^{2}+f^{2}(t)d\theta^{2}$ then
$H_{M}(s,\theta)=H_{\mathbb{M}(\kappa)}(s)$ for all $s\in (0,r]$
and all $\theta\in \mathbb{S}^{n-1}$ implies that $B_{M}(r)$ is
isometric to $B_{\mathbb{M}(\kappa)}(r)$. Because the equality
$H_{M}(s,\theta)=H_{\mathbb{M}(\kappa)}(s)$ for all $s\in (0,r]$
and all $\theta$ is equivalent  to have $\psi_{\kappa}(s)=0$,
$s\in [0,r]$ but this would imply that $f(s)=S_{\kappa}(s)$, $s\in
[0,r]$.
\end{remark}The next example shows that the rigidity in Theorem
(\ref{theoremCheng1}) is all we can expect without curvature
bounds.

\begin{example}\label{example1.1}For every $\kappa \in \mathbb{R}$, consider the  metric
$g=g(\kappa)$ on $M=[0, a]\times \mathbb{S}^{n-1}$, where
$a=\infty$ if $k\leq 0$ and $a =\pi/\sqrt{\kappa}$ if $\kappa>0$,
given in geodesic coordinates by the matrix  $g_{11}(t,\theta)=1$,
$g_{22}(t,\theta)= (S_{\kappa}^{4}(t)/t^{2})\cdot \theta_{22}$,
$g_{33}(t,\theta)=t^{2}\cdot\theta_{33} $, $g_{ii}(t,\theta)=
S_{\kappa}^{2}(t)\cdot \theta_{ii}$, $i\geq 4$,
$g_{ij}(t,\theta)=0$ if $i\neq j$, where
$d\theta^{2}_{ij}=(\theta_{ij})$ is the canonical metric of
$\mathbb{S}^{n-1}(1)$. This metric $g(\kappa)$ is smooth if
$\kappa\leq 0$. If $\kappa>0$ the metric $g(\kappa)$ is smooth
except at $(\pi,\theta)$. Let $h=h(\kappa)$ be the metric of
constant sectional curvature of $\mathbb{M}(\kappa)$ given by the
matrix $h_{11}=1$, $h_{ii}= S_{\kappa}^{2}(t)\cdot \theta_{ii}$,
 $i\geq2$. For $\kappa\neq 0$, $g(\kappa)$ is not isometric to
$h(\kappa)$.  Let $\triangle_{g}$ and $\triangle_{h}$ denote the
 Laplace operator of these two metrics  written in
 geodesic coordinates. They are given by
 \begin{eqnarray}\triangle_{g}&=&\frac{\partial^{2}}{\partial
 t^{2}}+ (n-1)\frac{C_{\kappa}}{S_{\kappa}}\frac{\partial}{\partial
 t}+\frac{t^{2}}{S_{\kappa}^{4}}\frac{\partial}{\partial
 \theta_{2}}+\frac{1}{t^{2}}\frac{\partial}{\partial
 \theta_{3}}+\sum_{i=4}^{n}\frac{1}{S_{\kappa}^{2}}\frac{\partial}{\partial
 \theta_{i}} \nonumber  \\
 &&\\
\triangle_{h}&=&\frac{\partial^{2}}{\partial
 t^{2}}+ (n-1)\frac{C_{\kappa}}{S_{\kappa}}\frac{\partial}{\partial
 t}+\sum_{i=2}^{n}\frac{1}{S_{\kappa}^{2}}\frac{\partial}{\partial
 \theta_{i}}.\nonumber
\end{eqnarray} We have that the
geodesic spheres $\partial B_{M}(s) $  and $\partial
B_{\mathbb{M}(\kappa)}(s)$ have the same mean curvature
$H_{M}(s)=H_{\mathbb{M}(\kappa)}(s)=(n-1)(C_{\kappa}/S_{\kappa})(s)$,
$s\in (0,r]$. And the geodesic balls $B_{M}(r) $ and $
B_{\mathbb{M}(\kappa)}(r)$ have the same first eigenvalue. For if
   $u$ be  a first Dirichlet eigenfunction of the geodesic ball
$B_{\mathbb{M}(\kappa)}(r)$, if $\kappa
>0$ suppose that that $r<\pi/\sqrt{\kappa}$. Thus
$\triangle_{h}u+\lambda_{1}(B_{\mathbb{M}(\kappa)}(r))u=0$ in $
B_{\mathbb{M}(\kappa)}(r)$  and $u=0$ on $\partial
B_{\mathbb{M}(\kappa)}(r)$. Since $u$ is radial we have that
 $\triangle_{h}u(t)=
\triangle_{g}u(t)=-\lambda_{1}(B_{\mathbb{M}(\kappa)}(r))u(t)$.
This shows that $u(t)$ is a first Dirichlet eigenfunction of the
geodesic ball $B_{M}(r)$  with same eigenvalue
$\lambda_{1}(B_{\mathbb{M}(\kappa)}(r)).$
\end{example}


\begin{thebibliography}{abcd}
\bibitem{kn:barta} Barta, J., {\em Sur la vibration fundamentale
d'une membrane.} C. R. Acad. Sci. 204,  472-473, (1937).

\bibitem{kn:bessa-montenegro2}Bessa, G. Pacelli and
Montenegro, J. F.,  {\em An Extension of Barta's Theorem \\
and geometric applications.} Deposited in ArchivX Paper Id:
math.DG/0308099, Preprint.

\bibitem{kn:chavel} Chavel, I., {\em Eigenvalues in Riemannian
Geometry.}  Pure and Applied Mathematics, (1984), Academic Press.

\bibitem{kn:cheng1} Cheng, S. Y., {\em Eigenfunctions and
eigenvalues of the Laplacian.} Am. Math. Soc. Proc. Sym. Pure
Math. 27, Part II, 185-193. (1975).

\bibitem{kn:cheng2} Cheng, S. Y., {\em Eigenvalue comparison
theorems and its geometric applications.} Math. Z. 143, 289-297,
(1975).



\bibitem{kn:veeravalli} Veeravalli, Alain R.:  {\em Une remarque
sur l'in\'{e}galit\'{e} de McKean.} Comment. Math. Helv. 78
884-888, (2003).


\end{thebibliography}
\end{document}